\documentstyle[12pt,amssymb,amstex]{amsart}
\newcommand{\cb}{{\cal B}}
\newcommand{\cc}{{\cal C}}
\newcommand{\cd}{{\cal D}}

\newcommand{\cf}{{\cal F}}

\newcommand{\ch}{{\cal H}}

\newcommand{\ck}{{\cal K}}
\newcommand{\cl}{{\cal L}}
\newcommand{\cam}{{\cal M}}

\newcommand{\co}{{\cal O}}
\newcommand{\cp}{{\cal P}}
\newcommand{\cq}{{\cal Q}}

\newcommand{\cx}{{\cal X}}

\newcommand{\ga}{{\frak A}}

\newcommand{\bc}{{\Bbb C}}

\newcommand{\bk}{{\Bbb K}}
\newcommand{\bm}{{\Bbb M}}
\newcommand{\bn}{{\Bbb N}}

\newcommand{\br}{{\Bbb R}}

\renewcommand{\a}{\alpha}
\renewcommand{\b}{\beta}
  
  \newcommand{\D}{\Delta}
\newcommand{\eps}{\varepsilon}

\newcommand{\r}{\rho}
\newcommand{\s}{\sigma}
\renewcommand{\t}{\tau}
\newcommand{\f}{\varphi} \newcommand{\F}{\Phi}
\newcommand{\th}{\theta}
\newcommand{\om}{\omega} \newcommand{\Om}{\Omega}

\newcommand{\itm}[1]{\item[($#1$)]}
\newcommand{\corr}[2]{{}_{#1}\hbox{id}_{#2}}
\newcommand{\id}{\hbox{id}}
\renewcommand{\ker}{\hbox{Ker}}
\newcommand{\re}{\hbox{Re}}
\newcommand{\tr}{\hbox{Tr}}

\newtheorem{thm}{Theorem}
\newtheorem{lem}{Lemma}
\newtheorem{cor}{Corollary}
\newtheorem{prop}{Proposition}

\theoremstyle{definition}
\newtheorem{defin}{Definition}

\theoremstyle{remark}
\newtheorem{rem}{Remark}
\newtheorem{step}{Step}

\begin{document}

\title[split property II]
{}
\begin{center}
{\bf OPERATOR SPACE STRUCTURES and the SPLIT PROPERTY II\\
The canonical non-commutative $L^2$ embedding} 
\end{center} 
\author{Francesco Fidaleo}
\address{Francesco Fidaleo\\
Dipartimento di Matematica\\
II Universit\`{a} di Roma (Tor Vergata)\\
Via della Ricerca Scientifica, 00133 Roma, Italy}
\email{{\tt Fidaleo@@mat.utovrm.it}}
\maketitle
 \begin{abstract}
A characterization of the split property for an inclusion $N\subset M$ of 
$W^*$-factors with separable predual is established in terms of the canonical
non-commutative $L^2$ embedding considered in \cite{B1,B2}
$$
\F_2:a\in N\to \D_{M,\Om}^{1/4}a\Om\in L^2(M,\Om)
$$
associated with an arbitrary fixed standard vector $\Om$ for $M$. 
This characterization 
follows an analogous characterization related to the canonical
non-commutative $L^1$ embedding
$$
\F_1:a\in N\to (\cdot\Om,J_{M,\Om}a\Om)\in L^1(M,\Om)
$$
also considered in \cite{B1,B2} and studied in \cite{F}. 
The split property for a Quantum Field Theory is 
characterized by equivalent conditions relative to the non-commutative
embeddings $\F_i$, $i=1,2$, constructed by the modular Hamiltonian of a 
privileged faithful state such as e.g. the vacuum state. The above 
characterization would be also useful for theories on a
curved space-time where there exists no a-priori privileged state.
 \end{abstract}

\section{Introduction}

The split property for inclusions of von Neumann algebras is of particular
interest from a theoretical viewpoint (\cite{B1,D,DL2}) although it was 
introduced and
intensively studied for the various applications in Quantum Field Theory,
see \cite{B0,B3,B5,Su}.
For the physical applications the relations of the split property
with several nuclearity conditions is also of main interest, 
see \cite{B2,B4,B6}.\\
The approach involving the study of properties of suitable sets
(or equally well suitable maps) was firstly introduced by Haag and
Swieca with the aim of characterizing those physical theories which are
asymptotically complete. They conjectured that
an asymptotically complete theory should satisfy a compactness criterion,
see \cite{HS}. Motivated by this approach, the split property 
was considered and a weak form of the celebrated Noether theorem
was established also for a quantized theory, see \cite{B3}. 
Moreover, for a theory with the split
property, the Haag-Swieca compactness criterion is automatically satisfied
\cite{B2}. A nuclearity conditition was proposed 
in \cite{B4}, as a stronger condition than the split property. 
A very relevant consequence of nuclearity is that theories satisfying
this property have a decent thermodynamical
behavior, see \cite{B5}.\\
To make the connection between the split property and the various
nuclearity conditions more transparent, canonical non-commutative embeddings
$M\hookrightarrow L^p(M)$, $p=1,2$ were considered firstly in \cite{B1}.
Given an inclusion $N\subset M$ of $W^*$-algebras with $M$ 
acting standardly on the
Hilbert space $\ch$ with a cyclic and separating vector $\Om$,
one can consider the following embeddings, canonically constructed
via the modular operators $\D$, $J$ associated to $\Om$, see \cite{B1}.
\begin{align*}
&\F_1:a\in N\to (\cdot\Om,Ja\Om)\in M_*,\\
&\F_2:a\in N\to \D^{1/4}a\Om\in \ch.
\end{align*}
In \cite{B1} it has been shown that, when $N\subset M$ is a factor-subfactor
inclusion, the nuclearity condition 
for any of the above maps insures the split property for the inclusion 
$N\subset M$. 
Conversely, if $N\subset M$ is a
split inclusion, these maps are nuclear for a dense set of cyclic separating
vectors for $M$, see \cite{B1}, Section 2. \\
Unfortunately the nuclearity condition is a stronger condition than
the split property if some privileged state (as e.g. the vacuum) is
kept fixed. Hence a complete characterization of the split property
in terms of properties of the $L^p$ embeddings would be desirable.
This approach has been followed in \cite{F} for the $L^1$ embedding 
where it has been shown that the metrically nuclear condition
for $\F_1$ characterize in a complete way the split property. In
this case canonical operator space structures on the involved 
normed spaces play a crucial role, see \cite{F}, Section 2.\\
The aim of this work is to establish a complete characterization
of the split property in terms of properties of the
non-commutative embedding $\F_2$. As an application one gets 
another characterization of the split property in Quantum
Field Theory even if the local algebras of observables have 
a non-trivial center, see e.g. \cite{B0,L}.
Therefore, in all the interesting cases arising from Quantum Field Theory,
the our characterization is available. To conclude this
introduction, we recall that the complete characterization
of the split property in terms of the above embeddings
$\F_i$, $i=1,2$ could be of interest also for theories living 
on curved space-time where there exists no a-priori privileged state as the
vacuum (\cite{W}). In this way the split property 
is directly stated in terms of properties
relative to the folium of states which are of interest for the theory.\\

This paper is organized as follows.\\

For the convenience of the reader we devote a preliminary part to resume 
some basic ideas on operator spaces considerated by several authors,
see \cite{B,E,E1,ER1,Pa,R}. A particular attention will be reserved 
to the structures recently
studied by Pisier such as the Hilbert $OH$ structure (\cite{P1,P3}) and the
non-commutative vector-valued $L^P$ spaces (\cite{P2,P4}). We also describe
the ideals $\cf_p(E,F)$, $p\geq1$, of the 
$p$-{\it factorable} maps between operator spaces $E$, $F$
considered in \cite{F1}. Successively, following \cite{B1}, we study
some basic properties of the above non-commutative $L^p$ embeddings
$\F_i:M\to L^p(M)$, $p=1,2$ relative to a $W^*$-algebra $M$.
As it has been well explained in \cite{B1} for a factor-subfactor inclusion
$N\subset M$, the extendibility of the
above embeddings, when restricted to the subfactor $N$,
is related to the split property for the given inclusion. 
We connect the extendibility property
of the above embeddings to a weaker property, the {\it quasi-split}
property, in the case of general inclusions
of $W^*$-algebras. The quasi-split property can also be stated using
the language of Connes correspondences \cite{Co1,CJ}.\\ 
A further section is devoted 
to study some extendibility properties of completely bounded normal maps of 
$W^*$-algebras with values
into a type $I$ factor as well as binormal
bilinear forms constructed in a natural way considering pairings
$$
a\otimes b\in M_1\otimes M_2\to\Psi_1\times\Psi_2(a\otimes b)
:=(\Psi_1(a),\overline{\Psi_2(b)})\in\bc
$$
where $\Psi_i:M_i\to\ch$ are normal maps of $W^*$-algebras with
values in a common Hilbert space $\ch$.
%Namely, if $\Psi_i$, $i=1,2$ are
%completely bounded (w.r.t. the Pisier $OH$ structure on $\ch$), 
%the functional $\Psi_1\times\Psi_2$ extends to a bounded binormal form
%on all of the projective $C^*$-tensor product $M_1\otimes_{max}M_2$.
These results enable us to prove the announced characterization
of the split property in terms of the canonical non-commutative
$L^2$ embedding $\F_2$:\\

\medskip
{\it Let $N\subset M$ be an inclusion of $W^*$-factors with separable 
predual and $\om\in M_*$ a faithful state.\\
The inclusion $N\subset M$ is a split inclusion if and only if the 
canonical $L^2$ embedding
$$
\F_{2,\om\lceil N}:a\in N\to\D^{1/4}a\Om\in L^2(M)
$$
is $2$-factorable as a map between the operator spaces
$N$ and $L^2(M)$ where $N$ has the natural operator space structure
as a $C^*$-algebra and $L^2(M)$ is equipped with the
Pisier structure $OH$.}\\

\medskip

As an immediate corollary
we have a geometrical characterization of the split property in terms
of the ``shape" of the image of the unit ball $\F_2(N_1)\subset L^2(M)$
under $\F_2$. Another
characterization of the split property suitable for the applications to
Quantum Field Theory is so stated.\\
We conclude with a section containing comments and problems concering
the natural embeddings
$\F_p:a\in A\to L^p(M)$
for a generic $p\geq1$ canonically associated to the inclusion $N\subset M$.\\
In this paper we deal only with $W^*$- algebras with separable
predual and all the operator spaces are considered to be complete if it is
not otherwise specified. For the general theory of operator algebras 
we refer the reader to the celebrated texts \cite{Sa,SZ,T}.

\section{Structures of operator spaces}

For the reader's convenience we collect 
some ideas about the operator
spaces which we need in the following. Details 
and proofs can be found later.\\

\subsection{Operator spaces}

We start with a normed space $E$, $E_1$ will be its (closed) unit ball.
Let $\|\ \|_n$ be a sequence of norms
on $\bm_n(E)$, the space of $n\times n$ matrices with entries in $E$. For
$a,b\in \bm_n$ these norms satisfy
\begin{align}
\begin{split}
&\|avb\|_n\leq\|a\|\ \|b\|\ \|v\|_n,\\
&\|v_1\oplus v_2\|_{n+m}=\max\{\|v_1\|_n,\|v_2\|_m\}
\label{ER}
\end{split}
\end{align}
where the above products are the usual row-column ones. 
This space with the above norms is called an (abstract){\it operator space}.
If\ $T:E\in F,\ T_n:\bm_n(E)\in \bm_n(F)$ are defined as
$T_n:=T\otimes\id$. $T$ is said to be {\it completely bounded} if 
$\sup\|T_n\|=\|T\|_{cb}<+\infty$; $\cam(E,F)$ denotes the set of all the 
completely bounded maps between $E$, $F$. It is an important 
fact (see \cite{R}) that a linear space $E$ with norms
on each $\bm_n(E)$ has a realization as a concrete operator space 
i.e. a subspace
of a $C^*$-algebra, if and only if these norms satisfy the properties 
in (\ref{ER}).\\
Given an operator space $E$ and $f=\bm_n(E^*)$, the norms
\begin{align}
\begin{split}
\|f\|_n=\sup&\{\|(f(v))_{(i,k)(j,l)}\|:v\in\bm_m(E)_1,\ m\in\bn\}\\
(&f(v))_{(i,k)(j,l)}=f_{ij}(v_{kl})\in \bm_{mn}
\label{dual}
\end{split}
\end{align}
determines an 
operator structure on $E^*$ that becomes itself an operator space (\cite{B}).
The linear space $\bm_I(E)$, $I$ any index set, is also of interest here.
$\bm_I(E)$ is in a natural manner an operator space via inclusion
$\bm_I(E)\subset\cb(\ch\otimes \ell^2(I))$ if E is realized as a subspace of
$\cb(\ch)$. Of interest is also the definition of $\bk_I(E)$ as those elements
$v\in\bm_I(E)$ such that $v=\lim\limits_\D v^\D$, where $\D$ is any finite 
truncation. Obviously
$\bm_I(\bc)\equiv\bm_I=\cb(\ell^2(I))$ and 
$\bk_I(\bc)\equiv\bk_I=\ck(\ell^2(I))$,
the set of all the compact operators on $\ell^2(I)$. For $E$ complete we remark
the bimodule property of $\bm_I(E)$ over $\bk_I$ because, for $\a\in\bk_I$,
$\a^\D v$, $v\a^\D$ are Cauchy nets in $\bm_I(E)$ and its limits
define unique elements $\a v$, $v\a$ that can be 
calculated via the usual
row-column product.\\
Given an index set $I$, we can define a map 
$\cx:\bm(V^*)\to\cam(V,\bm_I)$ which is a complete isometry, given by 
\begin{equation}
\label{mappa}
(\cx(f)(v))_{ij}:=f_{ij}(v).
\end{equation}
Moreover, if $f\in\bk(V^*)$, $\cx(f)$ is norm limit of finite rank
maps, so $\cx(\bk_I(V^*))\subset\ck(V,\bk_I)$, see \cite{ER2}, Section 3.\\
Remarkable operator space structures on a Hilbert space $H$ was introduced and
studied in \cite{ER3} via the following identification 
\begin{align*}
&\bm_{p,q}(H_c):=\cb(\bc^q,H^p),\\
&\bm_{p,q}(H_r):=\cb(\overline{H}^q,\bc^p)
\end{align*}
which define on $H$ the column and row structures $H_c$, $H_r$ respectively.\\
Recently the theory of complex interpolation has been developped
by Pisier for operator space as well. Successively a highly symmetric 
structure on a Hilbert space $H$ has been considered in \cite{P1,P3}, the
$OH$ structure. It can be viewed as an interpolating structure
between $H_c$, $H_r$:
$$
OH(I):=(H_c,H_r)_{1/2}
$$
where the cardinality of the index set $I$ is equal to the (Hilbert)
dimension of $H$. The explict descricption of the $OH$ structure, contained
in \cite{P3}, Theorem 1.1, is useful here and we report it in the next section 
for the sake of completness together with some others elementary properties.

\subsection{Tensor products between operator spaces}

Let $E$, $F$ be operator spaces, tensor norms on the algebraic tensor product 
$E\otimes F$ can be considerated. They make (the completions of) 
$E\otimes F$ operator spaces themselves.
The (operator) projective and spatial tensor products 
$E\otimes_{max}F$, $E\otimes_{min}F$ (\cite{ER1,ER2}) are of interest here
together with the Haagerup tensor product $E\otimes_h F$ (\cite{ER3}). 
As $F$ is complete, we have the completely isometric inclusion 
$E^*\otimes_{min}F\subset\cam(E,F)$. Also of interest is
the following complete identification 
\begin{equation}
\label{kappa}
\bk(E)=E\otimes_{min}\bk_I
\end{equation} 
for each operator space $E$.
The following important result concern the description of the predual of a 
$W^*$-tensor product in terms of the preduals of its individual factors.
\begin{thm}
let $N$, $M$ be any $W^*$-algebras. The predual
$(A\overline\otimes B)_*$ is completely isomorphic 
to the projective tensor product $A_*\otimes_{max}B_*$.
\end{thm}
The detailed proof of the above result can be found in \cite{ER2}, Section 3.
The above important result has a central role in the complete characterization
of the split property in terms of the $L^1$ embedding $\F_1$, see \cite{F}.

\subsection{The metrically nuclear maps}

In \cite{ER4} and, independently in \cite{F},
the class of the metrically nuclear maps $\cd(E,F)$ between operator
spaces $E$, $F$, has been introduced and studied. They can be
defined as
$$
\cd(E,F):=E^*\otimes_{max}F/\ker\cx
$$
where $\cx$ is the map (\ref{mappa}) which is 
a complete quotient map when restricted to $E^*\otimes_{max}F$. 
Also a geometrical characterization can be provided. 
Namely, an injective operator
$T\in\cd(E,F)$ is completely characterized by the shape of the
set $T(E_1)\subset F$, see \cite{F}, Section 2. The spaces $\cd(E,F)$ 
constituite, in the language of \cite{F1,Pi1}, an {\it operator ideal}
which have themselves a canonical operator space structure, see \cite{ER4,F}.

\subsection{The non-commutative vector-valued $L^p$ spaces}

The vector-valued non-commutative 
$L^p$ spaces has been introduced and intensively studied by Pisier 
\cite{P2,P4} by considering, for an operator space $E$, the
compatible couple of operator spaces
\begin{equation}
\label{comp}
\left(S_1(H)\otimes_{max}E,S_\infty(H)\otimes_{min}E\right);
\end{equation}
$S_1(H)$, $S_\infty(H)$ are the trace class and the
class of all compact operators acting on the Hilbert space $H$.
Then the vector-valued non-commutative $L^p$ spaces $S_p[H,E]$
can be defined as the interpolating spaces relative to the
compatible couple (\ref{comp})
\begin{equation}
\label{sp}
S_p[H,E]:=\left(S_\infty(H)\otimes_{min}E,S_1(H)\otimes_{max}E\right)_\th
\end{equation}
where $\th=1/p$.
As it has been described in \cite{P4} Theorem 1.1, the non-commutative
vector-valued $L^p$ spaces can be viewed as Haagerup tensor products
between the interpolating space structures on a Hilbert spaces
\begin{equation}
\label{sph}
S_p[H,E]=R(1-\th)\otimes_hE\otimes_h\overline {R(\th)}
\end{equation}
where $R(\th):=(H_r,H_c)_{\th}$ with $\th=1/p$.
In particular we have for $p=2$
\begin{equation}
\label{sp2}
S_2[H,E]=OH\otimes_hE\otimes_h\overline {OH}.
\end{equation}
The concrete structure of all the $S_p[H,E]$ can be derived from 
\cite{P2}, Theorem 2.
When there is no matter of confusion (i.e. if the Hilbert space is kept fixed),
we simply write $S_p[E]$ instead of $S_p[H,E]$.\\ 
For the reader's convenience we conclude with a result quite similar to 
that contained in \cite{ER2} Proposition 3.1 which will be useful
in the following, its proof can be found in \cite{F1}. 
We start with a Hilbert space
$H$ of (Hilbert) dimension given by the cardinality of the index set
$I$ and make the identification $H\equiv\ell^2(I)$.
\begin{prop}
\label{matr}
An element $u$ in $S_p[H,E]\subset\bm_I(E)$ satisfies $\|u\|_{S_p[H,E]}<1$
iff there exists elements $a,b\in S_{2p}(H)\subset\bm_I(E)$ with and
$v\in\bm_I(E)$ with $\|a\|_{S_{2p}(H)}=\|b\|_{S_{2p}(H)}=1$ and  
$\|v\|_{\bm_I(E)}<1$ such that
$$
u=avb.
$$
Furthermore one can choose $v\in\bk_I(E)$.
\end{prop}

\section{The Pisier $OH$ Hilbert space}

Here we collect some results concerning the Pisier $OH$ Hilbert space
which we need in the following.\\
We start with the concrete description of the $OH$ structure whose proof
can be easily recovered from \cite{P1}, Theorem 1.1.  
\begin{prop}
\label{oh}
Let $OH\equiv OH(I)$ be a Hilbert space equipped with the
$OH$ structure for a fixed index set $I$ and $x\in\bm_n(OH)$.
Then
$$
\|x\|_{\bm_n(OH)}=\|(x_{ij},x_{kl})\|^{1/2}_{\bm_{n^2}}
$$
where the entries of the above numerical matrix are as those in (\ref{dual}).
\end{prop}  
\begin{pf}
If $OH\subset\cb(\ch)$ then
$$
\|x\|_{\bm_n(OH)}=\|x\|_{\cb(\ch\otimes\bc^n)}.
$$
We consider $x$ in the dense set of $\bm_n(OH)$ consisting
of finite linear combinations of elements of the form
$x=a\otimes e_k$ where $a\in\bm_n$, $\{e_k\}_{k\in I}$
an ortonormal basis for $OH$. We have
$$
x_{ij}=\sum_{k\in I}(x_{ij},e_k)e_k
$$
where the coefficients in the above sum are all zero except
a finitely many of them. Then, by \cite{P3}, Theorem 1.1 part (iii),
we obtain
$$
\|x\|_{\bm_n(OH)}=\|\sum_{k\in I}a_k\otimes e_k\|_{\cb(\ch\otimes\bc^n)}
=\|\sum_{k\in I}
a_k\otimes{\overline a_k}\|^{1/2}_{\cb(\bc^n\otimes{\overline \bc^n})}
$$
where the numerical matrices $a_k$ are defined as 
$(a_k)_{ij}=(x_{ij},e_k)$. Therefore, as 
$(a_k\otimes{\overline a_k})_{(i,l)(j,m)}=(x_{ij},e_k)(e_k,x_{lm})$ we get
$$
\|x\|_{\bm_n(OH)}=\|\sum_{k\in I}(x_{ij},e_k)(e_k,x_{lm})\|^{1/2}
\equiv\|(x_{ij},x_{lm})\|^{1/2}_{\bm_{n^2}}.
$$
The proof now follows by a standard continuity argument.
\end{pf}
The following Proposition is a non-commutative version of the
Cauchy-Schwarz inequality and follows from a result due to
Haagerup.
\begin{prop}
\label{cs}
Let $OH\equiv OH(I)$ be as in the last Proposition
and $x\in\bm_m(OH)$, $y\in\bm_n(OH)$.
Then
$$
\|(x,y)\|_{\bm_{mn}}\leq
\|x\|_{\bm_m(OH)}\|\|y\|_{\bm_n(OH)}.
$$
\end{prop}  
\begin{pf}
We proceede as in the above Proposition.
$$
\|(x,y)\|_{\bm_{mn}}=
\|\sum_\s(x,e_\s)(e_\s,y)\|=
\|\sum_\s a_\s\otimes\bar b_\s\|_{\cb(\bc^m\otimes\bar \bc^n)}.
$$
Now we apply \cite{H}, Lemma 2.4 and obtain
\begin{align*}
&\|(x,y)\|_{\bm_{mn}}\equiv
\|\sum_\s a_\s\otimes\bar b_\s\|_{\cb(\bc^m\otimes\bar \bc^n)}\\
\leq&\|\sum_\s a_\s\otimes\bar a_\s\|^{1/2}_{\cb(\bc^m\otimes\bar \bc^m)}
\|\sum_\s b_\s\otimes\bar b_\s\|^{1/2}_{\cb(\bc^n\otimes\bar \bc^n)}\\
\equiv&\|x\|_{\bm_m(OH)}\|\|y\|_{\bm_n(OH)}
\end{align*}
which is the proof.
\end{pf}

\section{The class of factorable maps}

In this Section we resume the main properties of a class 
$\cf_p(E,F)\subset\ck(E,F)$, $1\leq p<+\infty$, of linear 
maps between operator spaces $E$, $F$ which are limits of finite
rank maps. These maps are obtained 
considering operators arising in a natural way from the Pisier 
non-commutative vector-valued $L^p$ spaces $S_p[H,E]$ and has been called
the {\it $p$-factorable} maps in \cite{F1}.
We also report for $p=2$, a geometrical description of the image
$T(E_1)\subset F$ of the unit ball of $E$ under an injective $2$-factorable
map. Details and proof can be found in \cite{F1}.

\subsection{The $p$-factorable maps}

\begin{defin}
Let $E$, $F$ be operator spaces and $1\leq p<+\infty$.
A linear map $T:E\in F$ will be called {\it $p$-factorable} if there
exists elements $b\in S_p[E^*]$, 
$A\in\cam\left(S_p,F\right)$ such that $T$ factorizes as
$$
T=AB
$$
where $B=\cx(b)\in\cam\left(E,S_p\right)$ and $\cx$ is the map 
(\ref{mappa}).\\
We also define
$$
\f_p(T):=\inf\{\|A\|_{cb}\|b\|_{S_p[E^*]}\}
$$
where the infimum is taken on all the factorization for $T$ as above.
The class of all the $p$-factorable maps between $E$, $F$ will be
denoted as $\cf_p(E,F)$.
\end{defin}

\medskip

In the above definition we have supposed $H\cong\ell^2$ 
without loss of generality, see \cite{F1}, Remark 1.
\begin{rem}
We always have $\cf_p(E,F)\subset\ck(E,F)$ as 
$\cx(\bk_I(E^*))\subset\ck(E,\bk_I)$.
\end{rem}
In \cite{F1} we have shown that 
$\left(\cf_p(E,F),\f_p\right)$, $1\leq p<+\infty$ are 
quasi-normed complete vector space. Moreover
we also have the ideal property for the above factorable maps, that is
$RST\in\cf_p(E_0,F_0)$ and
$$
\f_p(RST)\leq\|R\|_{cb}\f_p(S)\|T\|_{cb}
$$
whenever $E_0$, $E$, $F$, $F_0$ are operator spaces and $T:E_0\to E$,
$S:E\to F$, $R:F\to F_0$ linear maps with
$T\in\cam(E_0,E)$, $S\in\cf_p(E,F)$, $R\in\cam(F,F_0)$ respectively.\\
Actually the case with $p=2$ is particular: $\left(\cf_2(E,F),\f_2\right)$
is a Banach space for every operator space $E$, $F$.\\ 
Summarizing we have
\begin{thm}(\cite{F1} Theorem 2)
$(\cf_p,\f_p)$, $1\leq p<p$ are quasi-normed operator ideals
whereas $(\cf_2,\f_2)$ is a Banach operator ideal.
\end{thm}

\subsection{A geometrical description}

Analogously to the metrically nuclear operator setting \cite{F},
we give a suitable geometrical description for the range of a $2$-factorable
injective map. An analogous description can be stated also for a $p$-factorable
injective map, see \cite{F1}.\\
We start with an absolutely convex set $Q$ in an operator space
$E$ and indicate with $V$ its algebraic span. Consider a sequence
$\cq\equiv \{Q_n\}$ of sets such that

\medskip

\begin{itemize}
\item[(i)] $Q_1\equiv Q$ and each $Q_n$ is an absolutely convex
absorbing
set of $\bm_n(V)$ with $Q_n\subset\bm_n(Q)$;
\item[(ii)] $Q_{m+n}\cap(\bm_m(V)\oplus\bm_n(V))=Q_m\oplus Q_n$;
\item[(iii)] for $x\in Q_n$ then 
$x\in\lambda Q_n$ implies $bx\in\lambda Q_n, xb\in\lambda Q_n$
where $b\in(\bm_n)_1$.
\end{itemize}

\medskip

We say that a (possibly) infinite matrix $f$ with entries in the algebraic dual
of $V$ has finite $\cq$-norm if
$$
\|f\|_\cq\equiv\sup\{\|f^\Delta(q)
\|:q\in Q_n;\ n\in\bn;\ \Delta\}<+\infty
$$
where $f^\Delta$ indicates an arbitrary 
finite truncation corresponding to the finite set $\Delta$; 
the numerical matrix $f^{\Delta}(q)$ has entries as those in (1.2).

\medskip

\begin{defin}
\label{fatt}
An absolutely convex set $Q\subset E$ is said to be
{\it $(2,\cq)$-factorable} (where $\cq$ is a fixed 
sequence as above)
if there exists matrices $\a,\b\in S_4$ and a (possible infinite)
matrix $f$ of linear functionals as above with $\|f\|_\cq<+\infty$
such that, for each $x\in\cq_n$, one has
\begin{equation}
\label{ultra}
\|x\|_{\bm_n(E)}\leq C\|\a f(x)\b\|_{\bm_n(S_2)}.
\end{equation}
\end{defin}

\medskip

In the above case corresponding to $p=2$ we call a $(2,\cq)$-factorable
set simply $\cq$-factorable and omit the dependence on the sequence $\cq$
if there is no matter of confusion.\\
One can easily see that a $\cq$-factorable set is relatively compact,
hence bounded in the norm topology of $E$ and therefore $Q$, together the
Minkowski norms determinated by the $Q_n$'s on $\bm_n(V)$, is a (not necessarily
complete) operator space. Moreover the canonical immersion 
$V\stackrel{i}{\hookrightarrow}E$ is completely
bounded when $V$ is equipped with the operator structure determined
by the sequence $\cq$.\\

We now consider an injective completely bounded operator
$T:E\to F$ and the sequence $\cq_T$ given by $\cq_T=\{T_n(\bm_n(V)_1)\}$.
For such sequences the properties (i)--(iii) are automatically satisfied and, if
$T(V_1)$ is $\cq_T$-factorable, we call it simply $T$-factorable
and indicate the $\cq_T$-norm of a matrix of functional 
$f$ by $\|f\|_T$.
According to some interesting well-known cases, we have a description in
terms of the geometrical property of the range of certain $2$-factorable 
operators.
\begin{thm} 
Let $E$, $F$ operator spaces and $T:E\to F$ a
completely bounded injective map. $T\in\cf_2(E,F)$ iff $T(V_1)$ is a
$T$-factorable set in $F$.
\end{thm}
As in the case relative to the metrically nuclear maps,
the definition of a factorable set may appear rather involved; this is due
to the fact that the inclusion 
$\bm_n(V)_1\subset\bm_n(V_1)$ is strict in general
but, for an injective completely bounded 
operator $T$ as above, the $T$-factorable
set $T(V_1)$ is intrinsecally defined in terms of $T$.

\section{Inclusions of $W^*$-algebras and some related canonical embeddings}

We consider an inclusion $N\subset M$ of $W^*$-algebras with separable preduals
where $M$ acts standardly on $\ch$ with cyclic separating vector 
$\Omega$. Let
$\Delta,\ J$ be the Tomita's modular operator and conjugation relative to 
$\Omega$ respectively. Together
with $M$ we can consider its opposite algebra $M^\circ$ 
which is $*$-isomorphic to $M'$ via
$$
x\in M^\circ\in j(x)\equiv Jx^*J\in M';
$$
%\end{equation}
the inverse is given in the same way where 
$M^\circ$ is identified with $M$ as a linear space. The appearence of the 
opposite algebra $M^\circ$ will be clear later.\\
Following \cite{B1}, is of interest to consider the $L^1$ non-commutative 
embedding
\begin{equation}
\label{uno}
\Phi_1:b\in M\in(\cdot\Omega,Jb\Omega)\in L^1(M^\circ)
\end{equation}
where we identify in a natural manner $L^1(M^\circ)$ with $(M^\circ)_*$. 
The following non-commutative $L^2$ embedding is also of particular interest 
(\cite{B1,B2}).
\begin{equation}
\label{due}
\Phi_2:b\in M\in\Delta^{1/4}b\Omega\in L^2(M^\circ)
\end{equation}
where $L^2(M^\circ)\equiv \ch$.\\

If is not otherwise specified, all the Hilbert spaces appearing in the sequel
should be considered as endowed with 
the Pisier $OH$ structure as operator spaces.\\ 

It is a simple but important fact that the embeddings considerated above 
are completely positive in a sense which we are going to explain.\\

Let $A$, $B$ be $C^*$-algebras. A linear map $\Psi:A\to B^*$ is said to be
{\it completely positive} if all the maps
$$
\Psi\otimes\id:\bm_n(A)\to\bm_n(B^*)
$$
are positive that is $(\Psi\otimes\id)(a)$ is a positive element of 
$\bm_n(B)^*$ whenever $a$ is a positive element of $\bm_n(A)$.
It is easy to show that the complete positivity for $\Psi$
is equivalent to the request that the linear form $\om_\Psi$ given by
$$
\om_\Psi(a\otimes b):=\Psi(a)(b)
$$
define a positive form on the algebraic tensor product $A\otimes B$.\\

For the $L^2$ case we should start with a $W^*$-algebra and consider 
a hierarchy of self-polar cones as follows.
Let $M$ be a $W^*$-algebra with a normal semifinite faithful weight
$\f$. If $M$ is representated standarly on $L^2(M,\f)\equiv L^2(M)$
we can consider a hierarchy of self-polar cones $\cp_n\subset\bm_n(L^2(M))$,
$n\in\bn$, naturally associated to the normal semifinite faithful weight
$\f\otimes\tr_n$ on the $W^*$-algebras $\bm_n(M)$, see e.g. \cite{SZ} 10.23.
A characterization of when a hierarchy of self-polar cones 
$\cc_n\subset\bm_n(\ch)$, $n\in\bn$ is associated to a $W^*$-algebra
as above is given in \cite{S}.\\
Let $N$ be another $W^*$-algebra and $\F:N\to L^2(M)$ a linear map.
$\F$ will be said {\it completely positive} if the map 
$\F_n:\bm_n(N)\to \bm_n(L^2(M))$ is positive that is
$\F_n\left(\bm_n(N)_+\right)\subset\cp_n$ for each $n\in\bn$.\\

It is also of main interest to consider the normality property
of linear maps $T:M\to X$ of a $W^*$-algebra $M$ with values into 
an arbitrary linear space $X$.\\

Let $M$ be a $W^*$-algebra, $X$ a linear space with algebraic dual
$X'$ and $F\subset X'$ a
separating set of functionals. A continuous linear map $T:M\in X$ is
said to be {\it normal} ({\it singular}) w.r.t. $F$ 
according to the functionals $f\circ T$ on
$M$ are normal (singular) for every $f\in F$ (for the definition of 
normal and singular 
functionals on a von Neumann algebra see \cite{T}). Precisely, the $F$-normal
maps are just those which are $(\s(M,M_*),\s(X,F))$--continuous, whereas
the singular ones are those such that $T'(F)\subset M^\perp_*$ where $T'$
is the transpose map of $T$. From this
definition we get at once that a map which is normal and singular at the same
time must be the zero map, see \cite{T}. If $X$ is also a $W^*$-algebra,
it is usual to take $F\equiv M_*$; if $X$ is a normed space 
(which is not a $W^*$-algebra) and $F=X^*$, the topological dual, 
we call a normal map w.r.t. $X^*$ simply normal as well.\\

The normal property of the above non-commutative embeddings
are summarized in the following proposition (see \cite{B1}).
\begin{prop}
The non-commutative embeddings (\ref{uno}), (\ref{due}) 
are completely positive and normal.
\end{prop}
\begin{pf}
It is easy to show that $\F_i$, $i=1,2$ are completly positive 
and that $\F_1$ is normal so it remains to verify that $\F_2$ is 
normal too. Suppose that $\eta\in L^2(M^\circ)$, we obtain
\begin{align*}
|(\D^{1/4}&x\Om,\eta)|\leq\|\eta\|\|\D^{1/4}x\Om\|
=\|\eta\|(\D^{1/2}x\Om,x)^{1/2}\\
&\leq\|\eta\|(\|x\Om\|\|x^*\Om\|)^{1/2}
\leq{1\over\sqrt2}\|\eta\|(\|x\Om\|^2+\|x^*\Om\|^2)^{1/2}
\end{align*}
which shows that $\F_2$ is continuous in the strong$^*$ operator topology
hence $\s(M,M_*)$--continuous, see \cite{T}, Theorem II.2.6.
\end{pf}
A completely positive map between $C^*$-algebras is automatically
completely bounded, see e.g. \cite{Pa}; concerning the general case (i.e. when
the image space is not a $C^*$-algebra or merely an operator system)
it seems to be no relation between complete positivity and complete
boundedness. However the above embeddings are also completely
bounded.
\begin{prop}
The non-commutative embeddings (\ref{uno}), (\ref{due}) are completely bounded.
\end{prop}
\begin{pf} The first case is just \cite{F}, Proposition 3.1, so we prove 
the remaining one.
By Proposition \ref{oh}, we obtain for an element $a\in\bm_n(M)$
$$
(\D^{1/4}a_{ij}\Om,\D^{1/4}a_{kl}\Om)=(a_{ij}\Om,j(a^*_{kl})\Om)
=(a_{ij}j(a_{kl})\Om,\Om)\equiv\om_{n^2}(b)
$$
where $b\in\bm_{n^2}(\cb(L^2(M)))$ has norm equal to $\|a\|^2$.
Taking the supremum on the unit ball of $\bm_n(M)$ we obtain the assertion.
\end{pf}

\section{The quasi-split property for inclusions of $W^*$-algebras}

This Section follows Section 1 of \cite{B1} where the split property
for inclusions of factors is connected with the extendibility of
the canonical non-commutative embeddings $\F_i$, $i=1,2$
whose early properties have been summarized in the last Section.
Here we treat the general case of inclusions of $W^*$-algebra at all.
The results contained in this Section could be naturally applied
to inclusions arising from Quantum Field Theory where, under
general assumptions, the algebras of the local observables are
typically type $III$ algebras with (a-priori) a non-trivial 
center \cite{B0,L}.  This framework covers also
the commutative case which has been treated in some detail in
\cite{B1}, Section 3. We note that the subjects contained in this
Section have a natural description in terms of Connes {\it correspondence}
so we adopt somewhere the correspondence language in the following.\\ 
We start with an inclusion $N\subset M$ of $W^*$-algebras always
with separable predual, and suppose that $M$ is represented in standard
form on $L^2(M)$. We take a normal faithful state $\om\in M_*$ and
consider the Tomita operators $J,\D$ relative to it. 
\begin{defin}
The inclusion $N\subset M$ is said to be {\it quasi-split} if the map
\begin{equation}
\label{eta}
a\otimes b\in N\otimes M^\circ\to aJb^*J\in\cb(\ch)
\end{equation}
extends to a normal homomorphism $\eta$ of $N\overline{\otimes}M^\circ$
onto all of $N\vee M'$.
\end{defin}
We remark that, since the standard representation is 
unique up unitary equivalence,
the quasi-split property is really an intrinsic property of the inclusion.\\

Let $N$, $M$ be $W^*$-algebras, a $N-M$ {\it correspondence} is a separable 
Hilbert space $\ch$ which is a $N-M$-bimodule such that the 
(commuting) left and right 
actions of $N$, $M$ are normal.\\
For the basic facts
about the correspondence we remand the reader to \cite{Co1,CJ}.\\
The standard representation of $M$ gives rise to a correspondence which is
unique up unitary equivalence, just the identity $M-M$ correspondence
$\corr{M}{M}$ on $L^2(M)$ determined by the map (\ref{eta}).
If we have an inclusion $N\subset M$ of $W^*$-algebras we obtain a (uniquely
determined) $N-M$ correspondence if one restrict $\corr{M}{M}$
to $N$ on the left. We indicate this one as $\corr{N}{M}$ 
\begin{defin}
A $N-M$ correspondence $\s$ is said to be {\it split} if there
exists normal faithful representations $\pi$, $\pi^\circ$ of $N$, $M^\circ$
on Hilbert spaces $\ch_\pi$, $\ch_{\pi^\circ}$ respectively such that 
$\s$ is unitarily equivalent to $\pi\otimes\pi^\circ$
($\s\cong\pi\otimes\pi^\circ$ for short).
\end{defin}
The term {\it coarse} is used in \cite{Co1} for the case when 
$\pi$, $\pi^\circ$
are the standard representations of $N$, $M^\circ$ respectively.\\

Now we show that the quasi-split property can be viewed as a property of 
$\corr{N}{M}$.
\begin{prop}
Let $N\subset M$ an inclusion of $W^*$-algebras.\\
The following assertion are equivalent.
\itm{i} $N\subset M$ is a quasi-split inclusion.
\itm{ii} $\corr{N}{M}\prec\s$\\
where $\s$ is a $N-M$ split correspondence and $\prec$ means
the containment of representations.
\end{prop}
\begin{pf}
$(i)\Rightarrow(ii)$ If the map (\ref{eta}) extends to a normal 
homomorphism $\eta$ of $N\overline{\otimes}M^\circ$
onto $N\vee M'$ then $\eta=\eta_1\circ\eta_2\circ\eta_3$ where $\eta_3$
is an amplification, $\eta_2$ an induction and $\eta_1$ a spatial
isomorphism. Hence
$$
\eta(a\otimes b)=v^*(a\otimes b\otimes I)v
$$
where $v:L^2(M)\to L^2(M)\otimes L^2(M)\otimes\ch$ is an isometry that is
the $N-M$ correspondence $\corr{N}{M}$ is contained in a split one.\\
$(ii)\Rightarrow(i)$ If $\corr{N}{M}$ is a subcorrespondence
of a split one we have then a subrepresentation of the normal representation
$\pi\otimes\pi^\circ$ of $N\overline{\otimes}M^\circ$,
see \cite{T}, Theorem IV.5.2. Hence $\corr{N}{M}$ uniquely defines 
a normal representation of $N\overline{\otimes}M^\circ$ that
is a normal map $\eta$ of $N\overline{\otimes}M^\circ$ onto $N\vee M'$
which extends (\ref{eta}).
\end{pf}
When a faithful state $\om\in M_*$ is kept fixed, one can construct
the non-commutative embeddings $\F_1$, $\F_2$ given in 
(\ref{uno}), (\ref{due}) and, following \cite{B1},
one can study the extendibility properties of $\F_i$, $i=1,2$.\\

\vskip0.1cm

A completely positive normal map $\F:N\to L^p(M^\circ)$, $p=1,2$,
is said to be {\it extendible} (according with \cite{B1}, pag 236)
if, whenever $\widetilde{N}\supset N$ is another $W^*$-algebra
with separable predual, there exists a completely positive normal
map $\widetilde{\F}$ which yields commutative the following
diagram:

\begin{figure}[hbt]
\begin{center}
\begin{picture}(350,120)
\put(100,20){\begin{picture}(150,100)
\thicklines
\put(40,25){\vector(0,1){60}}
\put(47,91){\vector(1,-1){68}}
\put(50,15){\vector(1,0){62}}
\thinlines
\put(35,90){$\widetilde N$}
\put(35,10){$N$}
\put(118,10){$L^p(M^\circ)$}
\put(95,55){$\widetilde\F$}
\put(70,0){$\F$}
\end{picture}}
\end{picture}
\end{center}
\end{figure}

The following Theorem is the natural extension of \cite{B1}, 
Proposition 1.1.
\begin{thm}
Let $N\subset M$ be an inclusion of $W^*$-algebras and $\om\in M_*$
a fixed faithful state. Consider the embeddings $\F_i$, $i=1,2$ constructed
by $\om$ as in (\ref{uno}),(\ref{due}).\\
The following statements are equivalent.
\itm{i} $N\subset M$ is a quasi-split inclusion.
\itm{ii} $\F_{1\lceil N}:N\to L^1(M^\circ)$ is extendible.
\itm{iii} $\F_{2\lceil N}:N\to L^2(M^\circ)$ is extendible.
\end{thm}
\begin{pf}
$(iii)\Rightarrow(ii)$ We can consider the following commutative
diagram

\begin{figure}[hbt]
\begin{center}
\begin{picture}(350,120)
\put(100,20){\begin{picture}(150,100)
\thicklines
\put(33,75){\vector(1,0){84}}
\put(28,67){\vector(1,-1){36}}
\put(81,32){\vector(1,1){36}}
\thinlines
\put(20,70){$M$}
\put(53,20){$L^2(M^\circ)$}
\put(120,70){$L^1(M^\circ)$}
\put(70,83){$\F_1$}
\put(32,40){$\F_2$}
\put(110,40){$\Psi$}
\end{picture}}
\end{picture}
\end{center}
\end{figure}
where $\Psi$ given by
$$
\Psi(x):=(\D^{1/4}\cdot\Om,Jx)
$$
is computed by the transpose map of $\F_2$. So, if $\F_{2\lceil N}$
is extendible, it is easy to verify that $\F_{1\lceil N}$
is extendible too.\\
$(ii)\Rightarrow(i)$ Let $F$ be a type $I$ factor with separable predual
containing $N$ and $\widetilde{\F}_1$ the corresponding completely
positive normal extension of $\F_1$. Then
$$
\widetilde{\f}(f\otimes b):=\widetilde{\F}_1(f)(b)
$$
gives rise to a $F-M$ correspondence $(\s,\ch_\s)$.
Moreover the $F-M$ correspondence $(\s,\ch_\s)$, when restricted to
a $N$-$M$ correspondence, contains $\corr{N}{M}$. The assertion
now follows as the standard representation of $M$ is faithful.\\
$(i)\Rightarrow(iii)$ It is enough to show that $\F_2$ extends to a
type $I$ factor with separable predual containing $N$ 
(\cite{B1}, Proposition 1.1). 
Let $(\pi\otimes\pi^\circ,\ch_\pi\otimes\ch_{\pi^\circ})$ be a split
correspondence containing $\corr{N}{M}$. Then there exists
an isometry $v:L^2(M)\to\ch_\pi\otimes\ch_{\pi^\circ}$ with range projection
$vv^*\in \pi(N)'\overline{\otimes}\pi^\circ(M^\circ)'$, such that, for
$a\in N$, $b\in M^\circ$, we get
$$
aJb^*J=v^*\pi(a)\otimes\pi^\circ(b)v.
$$
As $vv^*\in\cb(\ch_\pi)\overline{\otimes}\pi^\circ(M^\circ)'$, it is easy to
show that $v^*(\cb(\ch_\pi)\otimes I)v\in M$. Now we define 
$\widetilde{\F}_2:\cb(\ch_\pi)\to L^2(M^\circ)$ as
$$
\widetilde{\F}_2(f):=\D^{1/4}v^*(f\otimes I)v\Om
$$
which is a completely positive normal extension of $\F_2$ to the type
$I$ factor $\cb(\ch_\pi)$.
\end{pf}
We note that the Proposition 1.1 of \cite{B1} can be immediately recovered
as a corollary of the above result.
\begin{cor}
If $N\subset M$ is an inclusion of $W^*$-factors, then the condition
$(i)$ in the last Theorem can be replaced by the a-priori stronger condition
\itm{i'} $N\subset M$ is a split inclusion.
\end{cor}
\begin{pf}
Let $\eta$ be the normal map which extends (\ref{eta}) to all of
$N\overline{\otimes}M$. If $N$, $M$, are both factors then 
$\ker\eta=\{0\}$ that is $\eta$ is a normal isomorphism
onto $N\vee M'$ but, in this case (\cite{D}, Corollary 1) that
condition turn out to be equivalent to the split property.
\end{pf}
As it has been explained in \cite{B1} for factor-subfactor inclusions,
the above theorem gets a characterization of the quasi-split property
in terms of the extendibility of the canonical embeddings
$\F_1$, $\F_2$. Moreover the extendibility
condition allows us to characterize the split property for
inclusions of factors (or also of algebras in some interesting cases such
those arising from Quantum Field Theory) directly in terms
of properties of $\F_1$, $\F_2$ somewhat similar to the nuclear
condition. This has been made in \cite{F} for the map $\F_1$.
We continue this program also for the $L^2$ embedding $\F_2$
which is also of interest for the applications, see \cite{B2}.

\section{On some completely bounded maps of 
$C^*$-algebras and $W^*$-algebras}

In this Section we analyze some extendibility properties of certain completely
bounded maps or bilinear (binormal) 
form on $C^*$ and $W^*$-algebras as well. These results will be crucial in the 
following. 
The results contained in this section are obtained under no separability 
conditions on the $W^*$-algebras which we are interested in.\\
We start with a result regarding the canonical structure 
of normal maps of a $W^*$-algebras with
values in an arbitrary type $I$ factor. This result directly follows from 
the analogous one contained in \cite{A} and relative to $C^*$-algbras.
\begin{prop}
\label{star}
Let $M$ be a $W^*$-algebra and $\ch$ a Hilbert space. Suppose that a 
completely bounded normal map $\pi:M\to\cb(\ch)$ is given.
\itm{i} There exists another Hilbert space $\ck$, bounded operators
$V_i:\ch\to\ck$, $i=1,2$ and a normal homomorphism $\r:M\to\cb(\ck)$
such that
$$
\pi=V_1^*\r(\cdot)V_2.
$$
\itm{ii} If $F\supset M$ is a type $I$ factor containing $M$, then $\pi$
extends to a completely bounded normal map $\widetilde\pi:F\to\cb(\ch)$.
\end{prop}
\begin{pf}
$(i)$ $\pi$ is completely bounded so there exists another 
Hilbert space $\ck$, bounded maps $V,W:\ch\to\ck$ 
and an homomorphism $\s:M\to\cb(\ck)$ such that 
$$
\pi=V^*\s(\cdot)W,
$$
see \cite{Pa},Theorem 7.4, or \cite{A} for the original exposition. Moreover we can decompose
$\s$ in its normal and singular part $\s=\s_n\oplus\s_s$ so we can write for
$V$, $W$
$$
V=
\begin{pmatrix}
V_n\\
V_s\\
\end{pmatrix},
\quad
W=
\begin{pmatrix}
W_n\\
W_s\\
\end{pmatrix}.
$$
We get for $\pi$
$$
\pi=W^*_n\s_n(\cdot)V_n+W^*_s\s_s(\cdot)V_s.
$$
As $\pi$ is normal, we have that 
$\pi-W^*_n\s_n(\cdot)V_n\equiv W^*_s\s_s(\cdot)V_s$ is, at the same time,
normal and singular and must be the zero map so the first part is proved
if one chooses $\r:=\s_n$.\\
$(ii)$ As $N\subset F\equiv\cb(\cl)$, the homomorphism
$\r$ given in $(i)$ can be written as 
$\r=\r_1\circ\r_2\circ\r_3$ where $\r_3$ is an
ampliation, $\r_2$ an induction and $\r_1$ a spatial isomorphism.
Namely there exists a suitable Hilbert space $\widehat\cl$ and 
an isometry $U:\ck\to\cl\otimes\widehat\cl$ such that
$$
\r(a)=U^*(a\otimes I)U.
$$
Now, if $f\in\cb(\cl)$, we define
$$
\widetilde\pi(f):=V^*_1U^*(f\otimes I)UV_2
$$
which is a completely bounded normal extension of $\pi$ to all
of $F$.
\end{pf}
Another interesting situation is when one consider a pairing constucted 
in a natural manner by linear maps of $C^*$-algebras or $W^*$-algebras
with values in a Hilbert space equipped with the Pisier $OH$ structure
as an operator space.\\
Let $M_i$, $i=1,2$ be $C^*$-algebras and $\Psi_i$, $i=1,2$ linear
maps with range in a fixed Hilbert space $\ch$. Then a
bilinear form $\Psi_1\times\Psi_2:M_1\otimes M_2\to\bc$ is defined as follows
\begin{equation}
\label{sette}
(\Psi_1\times\Psi_2)(a_1\otimes a_2):=(\Psi_1(a_1),\overline{\Psi_2(a_2)})
\end{equation}
where we have made the usual (antilinear) identification between $\ch$
and $\overline\ch$. Suppose that $\ch$ is endowed with the Pisier $OH$ structure
as an operator space. We have the following
\begin{prop}
\label{lim}
If $\Psi_i$, $i=1,2$ are completely bounded then the linear form (\ref{sette})
extends to a bounded functional on the projective tensor product
$M_1\otimes_{max}M_2$.\\ 
Moreover if $M_i$, $i=1,2$ are 
$W^*$-algebras and $\Psi_i$, $i=1,2$ are normal, then the
linear form (\ref{sette}) extends to a binormal bounded form on all of
$M_1\otimes_{max}M_2$.
\end{prop}
\begin{pf}
By \cite{ER1} Theorem 3.2, the above bilinear form extends
to a bounded functional on $M_1\otimes_{max}M_2$ iff 
$\|\Psi_1\times\Psi_2\|_{cb}<+\infty$. We get by Proposition \ref{cs}
$$
\|(\Psi_1(a),\overline{\Psi_2(b)}\|_{\bm_{mn}}\leq
\|\Psi_1(a)\|_{\bm_m(OH)}\|\Psi_2(b)\|_{\bm_n(OH)}\leq
\|\Psi_1\|_{cb}\|\Psi_2\|_{cb}
$$
and, taking the supremum on the left on the unit balls of $M_1$, $M_2$,
we obtain the first part. The last part directly follows by definition.
\end{pf}

\section{The split property for inclusions of $W^*$-Algebras and 
canonical non-commutative embeddings}

In this section we provide the announced complete characterization
of the split property for a factor--subfactor inclusion $N\subset M$
in terms of certain factorability property of the non-commutative
$L^2$ embedding ${\displaystyle \F_{2\lceil N}:N\to L^2(M)}$ given
in (\ref{due}). The final part of this Section is devoted to extend
this characterization to inclusions of $W^*$-algebras which arise
from Quantum Field Theory and have a-priori a non-trivial center,
see e.g. \cite{B0,L}.\\

The following Lemma is a starting point.
\begin{lem}
\label{tre}
Let $M$ be a $W^*$-algebra and $V$ an operator
space. If $T\in\cf_p(M,V)$ is a normal map then $T$ has a decomposition
$$
T=A\cx(b)
$$
where $b$ can be chosen in $S_p[M_*]$ and $\|A\|_{cb}\|b\|_{S_p[M^*]}$
is arbitrary close to $\f_p(T)$.
\end{lem}
\begin{pf} 
Suppose that $Tx=A\a f(x)\b$ with $f\in\bm_\infty(M^*)$ with
$\|f\|\leq\f_p(T)+\eps$. We can decompose $f$ in its normal and
singular parts $f^n,\ f^s$; then we have
$$
Tx=A\a f^n(x)\b+A\a f^s(x)\b
$$
where the summands on the l.h.s. are well defined, 
respectively normal and singular,
operators in $\cf_p(M,V)$ as $\|f^n\|,\|f^s\|\leq\|f\|$. We than have that 
$T-A\a f^n(\cdot)\b$
is, at the same time, normal and 
singular and must be the zero mapping.
Therefore $$
T=A\a f^n(x)\b
$$
and $\f_p(T)\leq\|f^n\|\leq\|f\|\leq\f_p(T)+\eps$.
\end{pf}
Now we are ready to prove the characterization of the split
property for a factor--subfactor inclusion.
\begin{thm}
\label{mega}
Let $N\subset M$ be an inclusion of $W^*$-factors with 
separable preduals and $\om\in M_*$ a faithful state. Let
$\F_i:M\to L^i(M^\circ)$, $i=1,2$ be the embeddings associated to a normal
faithful state $\omega$ for $M$ and given in (\ref{uno}), (\ref{due}).\\ 
The following statements are equivalent.
\itm{i} $N\subset M$ is a split inclusion.
\itm{ii} $\F_{1\lceil N}\in\cd(N,(L^1(M^\circ))$.
\itm{ii'} The set $\{(\cdot\Omega,Ja\Omega):a\in A,\ \|a\|<1\}$ is
$\F_1$-decomposable (see \cite{F}, Definition 2.6 for this definition).
\itm{iii} $\F_{2\lceil N}\in\cf_2(N,(L^2(M))$.
\itm{iii'} The set $\{\D^{1/4}a\Om:a\in A,\ \|a\|<1\}$ is
$\F_2$-factorable.
\end{thm}
\begin{pf}
Some of the above equivalences are immediate (\cite{F1})
or are contained in \cite{F} so we prove the remaining ones.\\
$(i)\Rightarrow(iii)$ If there exists a type $I$ interpolating factor $F$
then $\F_2$ factors according to

\newpage

\begin{figure}[hbt]
\begin{center}
\begin{picture}(350,120)
\put(100,20){\begin{picture}(150,100)
\thicklines
\put(33,75){\vector(1,0){84}}
\put(28,67){\vector(1,-1){36}}
\put(81,32){\vector(1,1){36}}
\thinlines
\put(20,70){$N$}
\put(60,20){$L^2(F)$}
\put(120,70){$L^2(M)$}
\put(70,83){$\F_2$}
\put(32,40){$\Psi_2$}
\put(115,40){$\Psi_1$}
\end{picture}}
\end{picture}
\end{center}
\end{figure}

where $\Psi_2$ arises from $S_2[N_*]$ and $\Psi_1$ is bounded, see \cite{B1}. 
Moreover $\Psi_1$ is automatically completely bounded, see \cite{P3}, 
Proposition 1.5.\\
$(iii)\Rightarrow(i)$ It is enough to show that, if $N\subset F$ with
$F$ a type $I$ factor with separable predual, 
then $\F_1$ extends to a completely positive
map $\widetilde{\F}_1:F\to L^1(M^\circ)$, see \cite{B1}, Proposition 1.1.
For the reader's convenience we split up this part of the proof in several
steps.\\

\begin{step}
If ${\displaystyle \F_{2\lceil N}\in\cf_2(N,M)}$ then 
${\displaystyle\F_{2\lceil N}=A\cx(b)}$
where $b$ can be chosen in $S_2[N_*]$.
\end{step}
\begin{pf}
This directly follows by Lemma \ref{tre} as $\F_2$ is a normal map.
\end{pf}
\begin{step}
If $F\supset N$ is a Type $I$ factor with separable predual then
$\F_2$ extends to a completely bounded normal map
$\widetilde{\F}:F\to L^2(M^\circ)$.
\end{step}
\begin{pf}
As $\F_2(x)=A\a f(x)\b$ for $x\in N$, where $f\in\bm_\infty(N_*)$,
we get a completely bounded normal map $\r:=\cx(f):N\to\bm_\infty$.
By Proposition \ref{star}, if $F\supset N$ is a Type $I$ factor with 
separable predual we have a completely bounded normal extension
$\tilde\r:F\to\bm_\infty$ to all of $F$. If we take 
$\tilde f:=\cx^{-1}(\tilde\r)$ we obtain the desired extension.
\end{pf}
Unfortunately $\widetilde{\F}$ might be not positive, so we should
recover a completely positive normal extension of $\F_{1\lceil N}$
to $F$ via $\widetilde\F$.

\begin{step}
Let $\f:F\otimes M^\circ\to\bc$ be the binormal bilinear form given by
$$
\f(x\otimes y):=(\widetilde{\F}(x),\D^{1/4}y^*\Om).
$$
Then $\f$ uniquely defines a bounded binormal form on all the $C^*$-algebra
$F\otimes_{max}M^\circ$.
\end{step}
\begin{pf}
As $y\in M^\circ\to y^*\in M$ is a completely bounded antilinear map,
we obtain the assertion by Proposition \ref{lim}.
\end{pf}
\begin{step}
$\f$ can be decomposed in four positive binormal form.
\end{step}
\begin{pf}
We consider the universal envelopping von Neumann algebra $\cam$ of 
$F\otimes_{max}M^\circ$ together with the central projection $p$
relative to the binormal forms as described in \cite{T}, Theorem
III.2.7. Then $\f$ gets a normal form on the von Neumann algebra
$\cam_p$ and we can recover the desired decomposition via the 
Jordan decomposition described in \cite{T}, Theorem III.4.2.
\end{pf}

Let $\f_+$ be the positive part of $\f$, then $\f_+$ gives rise to a
cyclic $F$--$M$ correspondence.\\

\begin{step}
$\f_+$, when restricted to $N\otimes M^\circ$, dominates $\om$ given by
$$
\om(x\otimes y):=(x\Om,Jy\Om).
$$
\end{step}
\begin{pf}
As $\re\f=\f_+-\f_-$ we get
$$
\f_+(a)=\om(a)+\f_-(a)\geq\om(a)
$$
whenever $a\in N\otimes M^\circ$.
\end{pf}

Now we consider the $GNS$ construction for the cyclic $F-M$ correspondence
determined by $\f_+$ and obtain two
normal commuting representations $\pi$, $\pi^\circ$ of $F$, $M^\circ$
respectively on a separable Hilbert space $\ch$ and a vector
$\xi\in\ch$ (cyclic for $\pi(F)\vee\pi^\circ(M^\circ)$) such that
$$
\f_+(x\otimes y)=(\pi(x)\pi^\circ(y)\xi,\xi).
$$
Then, as the restriction of $\f_+$ dominates $\om$, 
we have a positive element $T\in\pi(N)'\wedge\pi^\circ(M^\circ)'$
such that, if $x\in N$, $y\in M^\circ$
$$
\om(x\otimes y)=(\pi(x)\pi^\circ(y)T\xi,T\xi).
$$
Now we define $\widetilde{\F}_1:F\to L^1(M^\circ)$ given by
$$
\widetilde{\F}_1(f):=(T\pi(f)T\pi^\circ(\cdot)\xi,\xi)
$$
which is a completely positive normal map which extends
$\F_1$. This completes the proof.
\end{pf}
As there exists examples of quasi-split inclusions which do not satisfy
the split property (\cite{D}), it is still unclear if the former can
be characterized via the $2$-factorable maps also for the
case of inclusions of $W^*$-algebras with a non-trivial center. However in some
interesting cases such as those arising from Quantum Field Theory
we turn out to have the same characterization.
We suppose that the net $\co\to\ga(\co)$ of von Neumann algebras of local
observables of a quantum theory acts on the Hilbert space $\ch$ 
and satisfy all the usual assumptions (a priori
without the split property) which are typical in Quantum Field Theory,
$\Om\in\ch$ will be the vacuum vector which is cyclic for the net, 
see e.g. \cite{Su}. We prefer to characterize the split property
directly in terms of the shape of the sets $\F_i(\ga(\co)_1)$, $i=1,2$
which are of particular interest in Quantum Field Theory, 
see \cite{B2,B4,B5,B6,HS}.
\begin{thm}
Let $\co\subset int(\widehat\co)$ be double cones in the physical
space-time and $\ga(\co)\subset\ga(\widehat\co)$ the corresponding
inclusion of von Neumann algebras of observables.\\
The following assertions are equivalent.
\itm{i} $\ga(\co)\subset\ga(\widehat\co)$ is a split inclusion.
\itm{ii} The set $\{(\cdot a\Om,\Om):a\in\ga(\co)_1\}
\subset(\ga(\widehat\co)')_*$ is
$\F_1$-decomposable.
\itm{iii} The set $\{\D^{1/4}a\Om:a\in\ga(\co)_1\}
\subset\ch$ is $\F_2$-factorable.
\end{thm}
\begin{pf}
If $\ga(\co)\subset\ga(\widehat\co)$ is a split inclusion then $(ii)$
and $(iii)$ are true, see \cite{B1,F,F1}. Conversely, if $(ii)$
or $(iii)$ are satisfied, then $\F_1$ is extendible, see \cite{F} for the
implication $(ii)\Rightarrow (i)$ or the proof of the last Theorem 
for the implication $(iii)\Rightarrow (i)$. Hence the map 
$$
\eta:a\otimes b\in\ga(\co)\otimes\ga(\widehat\co)'\to 
ab\in\ga(\co)\vee\ga(\widehat\co)'
$$
extends to a normal homomorphism of 
$\ga(\co)\overline{\otimes}\ga(\widehat\co)$ onto 
$\ga(\co)\vee\ga(\widehat\co)$. But this homomorphism 
is in fact an isomorphism by an argument exposed in \cite{B0}, pagg. 129-130.
Moreover, as $\ga(\co)\wedge\ga(\widehat\co)$ is properly infinite
(\cite{Ka}), the assertion now follows by Corollary 1 of \cite{D}.
\end{pf}
For the applications to Quantum Field Theory, the split property turns
out to be completely characterized by properties regarding the
non-commutative embeddings $\F_{i,\om}$, $i=1,2$ given in (\ref{uno}),
(\ref{due}) which do not depend on the choice of the faithful state
$\om\in M$. We wish to note that the characterizations listed above 
could be of interest in particular for the applications to
theories living on a curved space-time where there is no privileged
state as the vacuum, see e.g. \cite{W}.

\section{Remarks about the canonical embeddings $M\hookrightarrow L^p(M)$ 
with $p$ arbitrary}

In this paper we have studied properties relative to symmetric
embedding $\F_2:M\to L^2(M)$. Moreover a complete characterization
of the split property has been given for an inclusion $N\subset M$. 
If one collect the above results
together with those contained in \cite{B1,F}, we have a situation
relative to the embeddings $\F_i$, $i=1,2$ which seems to be satisfactory 
enough.
On the other hand, it would be of interest to consider other canonical
embeddings which can be constructed for arbitrary $p$. 
This can be made in the following way. We leave fixed
a faithful state $\om\in M_*$ and consider the crossed product
$R:=M\times_{\s^\om}\br$ where $\s^\om$ is the modular group relative to
$\om$. $R$ has a faithful semifinite normal trace $\t$ and the dual 
action $\th$ (scaling automorphisms) of $\br$ satisfies 
$\displaystyle{\t\circ\th_s=e^{-s}\t}$. If $1\leq p<+\infty$,
one can identify the $L^p$ spaces with the spaces of all $\t$-measurable
operators $h$ affiliated to $R$ which satisfy, for $s\in\br$
$$
\th_s(h)=e^{-{s\over p}}h.
$$
If we fix a vector $h_\om\in L^2$ which represents $\om$, the canonical 
symmetric non-commutative embeddings
$\F_p:M\to L^p(M)$, $1\leq p<+\infty$ assume the following simple form 
\begin{equation}
\label{generico}
\F_p(x):=h_\om^{{1\over 2p}}xh_\om^{{1\over 2p}}.
\end{equation}
It easy to show that the last definition coincide with our previous one
in the case when $p=1,2$. As a first step, we have shown that, if $p=1,2$,
the above embeddings are completely bounded when a suitable
operator space structure is kept fixed on $L^p(M)$. These are precisely

\medskip

\begin{itemize}
\item[(a)] the canonical (pre)dual structure of $L^1(M)$ 
when it is considerated as the 
predual of $M^\circ$, the opposite algebra of $M$, for the embedding $\F_1$,
\item[(b)] the Pisier $OH$ structure on $L^2(M)$ for the embedding
$\F_2$.
\end{itemize}

\medskip

Moreover, properties of restrictions $\F_{1\lceil N}$, $\F_{2\lceil N}$,
characterize the split property for the inclusion $N\subset M$. Also
in this case the operator space structures considerated on $L^1(M)$, $L^2(M)$
play a crucial role. On all of $L^p(M)$, it will be a lot of
interpolating operator space structures. Unfortunately this picture is
not yet understood in the full generality, see \cite{P4}. In the light
of the above considerations, it would be of interest to extend 
our characterization of the split property in terms of the other
embeddings $\F_p:M\hookrightarrow L^p(M)$ (\ref{generico}). We remark that also
non-symmetric embeddings could be considerated as well, see \cite{K1}, 
Section 7. We hope to return on these interesting questions somewhere else.

\vspace{1.5cm}

\end{document}